\documentclass[12pt]{article} 
\usepackage{amsmath}
\usepackage{amssymb} 
\usepackage{amscd}
\usepackage{amsthm}

\usepackage[utf8]{inputenc}
\usepackage{pb-diagram}
\usepackage{hyperref}

\newtheorem{theorem}{Theorem}

\newtheorem{proposition}[theorem]{Proposition}
\newtheorem{lemma}[theorem]{Lemma}

\def\ad{\mathrm{ad}}

\def\im{{\mbox{Im}}}
\def\dim{{\mbox{dim}}}

\def\Tr{{\mbox{Tr}}}

\def\ker{{\mbox{Ker}}}

\def\gldim{{\mbox{gldim}}}

\def\Ext {{\mbox{Ext}}}

\def\End{{\mbox{End}}}

\def\cala{{\mathcal A}} 
\def\calb{{\mathcal B}} 
\def\calk{{\mathcal K}}

\def\calt{{\mathcal T}}
\def\calc{{\mathcal C}}

\def\pfraca{{\mathfrak a}}
\def\pfracb{{\mathfrak b}}

 \def\fracg{{\mathfrak g}}

 \def\fracsl{{\mathfrak {sl}}}
\def\fracso{{\mathfrak {so}}}

\def\bbbone{\mbox{\rm 1\hspace {-.6em} l}}

\def\Tor{{\mbox{Tor}}}

\numberwithin{equation}{section}

\begin{document}

\enlargethispage{3cm}

\thispagestyle{empty}
\begin{center}
{\bf QUADRATIC  ALGEBRAS ASSOCIATED}\\
\vspace{0.2cm}
{\bf  WITH EXTERIOR 3-FORMS}\\
\end{center}
   
\vspace{0.3cm}

\begin{center}
Michel DUBOIS-VIOLETTE
\footnote{Laboratoire de Physique des 2 Infinis Ir\`ene Joliot Curie\\
	P\^ole Th\'eorie, 
	IJCLab UMR 9012\\
	CNRS, Universit\'e Paris-Saclay, 
	B\^atiment 210\\
	F-91406 Orsay Cedex, France\\
	michel.dubois-violette@u-psud.fr} and 
Blas TORRECILLAS
\footnote {Universidad de Almeria\\ Departemento de Matem\'aticas\\ Carretera Sacramento s/n\\ 04120 La Can\~ada de San Urbano, Almeria, Espa\~na\\
blas.torrecillas@gmail.com} \end{center}
\vspace{0,5cm}

\begin{abstract}
This paper is devoted to the study of the quadratic algebras with relations generated by superpotentials which are exterior  3-forms. Such an algebra is regular if and only if it is Koszul and is then a 3-Calabi-Yau domain. After some general results we investigate the case of the algebras generated in low dimensions $n$ with $n\leq 7$. We show that whenever the ground field is algebraically closed all these algebras associated with 3-regular exterior 3-forms are regular and are thus 3-Calabi-Yau domains. This result does not generalize to dimensions $n$ with $n\geq 8$  : we describe a counter example in dimension $n=8$.

\end{abstract}

\vfill
\today

 \newpage
\tableofcontents

\newpage
\section{Introduction}

Throughout $\mathbb K$ is a field and all vector spaces, algebras, tensor products, etc. are over $\mathbb K$. By an algebra without other specifications we mean a unital associative $\mathbb K$-algebra. We use the Einstein convention of summation over repeated up-down indices in the formulas.\\

Let $\cala=\oplus_{n\geq0}$ $\cala_n$ be a $\mathbb N$-graded connected algebra, $\cala_0=\mathbb K\bbbone$ is the trivial module and is identified with $\mathbb K$. The projective dimension of the trivial module is the global dimension of $\cala$ and is also its Hochschild dimension (in homology as well as in cohomology),\cite{car:1958},  \cite{ber:2005}. The algebra $\cala$ is said to be {\sl regular} if its global dimension is finite say $\gldim(\cala)=D$ $(\in \mathbb N)$ and if
\begin{equation}
\Ext^k(\mathbb K,\cala)=\left\{
\begin{array}{ccc}
\mathbb K& \text{if}& k=D\\
0 & \text{if} &  k\not= D
\end{array}
\right.
\end{equation}
i.e. $\Ext^k(\mathbb K,\cala)=\delta^k_D\mathbb K$.\\

In the following we shall be concerned about the connected algebras freely finitely generated in degree 1 with a finite number of relations $r_i$ of degrees $\geq 2$, thus by algebras of the form
\begin{equation}
\cala=\mathbb K\langle x^1,\dots, x^n\rangle / [\{r_i\}]
\end{equation}
where $\mathbb K\langle x^1,\dots, x^n\rangle$ is the free connected algebra generated by the $x^k$ $(k\in \{1,\dots, n\})$ and where $[S]$ denotes the ideal of $\mathbb K\langle x^1,\dots, x^n\rangle$  generated by $S\subset \mathbb K\langle x^1,\dots, x^n\rangle$. 
Notice that $\mathbb K\langle x^1,\dots, x^n\rangle$ is canonically isomorphic to the tensor algebra $T(\mathbb K^n)$ of $\mathbb K^n$. When all the relations $r_i$ are of the same degree $N$, $\cala$ is said to be an {\sl $N$-homogeneous algebra}. One has the following result \cite{ber-mar:2006}.

\begin{proposition}\label{reg}
Let $\cala$ be a regular algebra of global dimension $D$.\\
$(\imath)$ If $D=2$ then $\cala$ is quadratic and Koszul.\\
$(\imath\imath)$ If $D=3$ then $\cala$ is $N$-homogeneous with $N\geq 2$ and Koszul.
\end{proposition}
The notion of Kozulity introduced in \cite{pri:1970} for quadratic algebras has been extended in \cite{ber:2001a} for $N$-homogeneous algebras with $N\geq 2$.

For the case $D=2$ one has the following complete description \cite{zha:1998}.
\begin{proposition}\label{D2}
Let $\cala$ be a regular algebra of global dimension 2 then
\begin{equation}
\cala=\mathbb K\langle x^1,\dots, x^n\rangle/[b_{ij}x^i\otimes x^j]\label{3}
\end{equation}
where $b_{ij} =b(e_i,e_j)$ are the components of a nondegenerate bilinear  form $b$ on $\mathbb K^n$. Conversely if $b$ is a nondegenerate bilinear form on $\mathbb K^n$, Formula {\rm(\ref{3})} defines a regular algebra of global dimension 2.
\end{proposition}

In order to state the similar result of the first part of the above proposition for the case of the global dimension $D=3$, let us remind some definitions of \cite{mdv:2007} concerning multilinear forms on $\mathbb K^n$.\\

Let $Q\in GL(n,\mathbb K)$ and $m$ be an integer with $n\geq m\geq 2$; Then an $m$-linear form $w$ on $\mathbb K^n$ is said to be $Q$-{\sl cyclic} if one has
\begin{equation}
w(X_1,\dots,X_m)=w(QX_m,X_1,\dots, X_{m-1})
\end{equation} 
for any $X_1,\dots, X_m\in \mathbb K^n$. An $m$-linear form $w$ on $\mathbb K^n$ is said to be {\sl preregular} if it satisfies the following conditions $(\imath)$ and $(\imath\imath)$ :\\
$(\imath)$ $w(X,X_1,\dots, X_{m-1})=0$ for any $X_1,\dots, X_{m-1}\in \mathbb K^n$ implies $X=0$,\\
$(\imath\imath)$ there is an element $Q_w\in GL(n,\mathbb K)$ such that $w$ is $Q_w$-cyclic.\\

In view of $(\imath)$, $Q_w$ is then unique and this twisted cyclicity implies that whenever $w(X_1,\dots,X_k,X,X_{k+1},\dots, X_{m-1})=0$ for any $X_1,\dots, X_{m-1}\in \mathbb K^n$ then $X=0$; this latter condition will be referred to as 1-{\sl site nondegeneracy}.\\

It is worth noticing here that a bilinear form on $\mathbb K^n$ is preregular if and only if it is nondegenerate. Finally we need another definition of \cite{mdv:2007}, namely an $m$-linear form $w$ on $\mathbb K^n$ will be said to be 3-{\sl regular} if it is preregular and satisfies\\
$(\imath\imath\imath)$ If $L_1, L_2\in \End(\mathbb K^n)$ are such that
\[
w(L_1X_1,X_2,X_3,\dots,X_m)=w(X_1,L_2X_2,X_3,\dots,X_m)
\]
for any $X_1,\dots,X_m\in \mathbb K^n$, then $L_1=L_2=\lambda\bbbone$ for some $\lambda\in \mathbb K$.\\

 We can now formulate the global dimension 3 version of the first part of Proposition \ref{D2}, \cite{mdv:2005},\cite{mdv:2007}, \cite{mdv:2010}.

\begin{proposition}\label{D3}
Let $\cala$ be a regular algebra of global dimension 3 then
\begin{equation}
\cala=\mathbb K\langle x^1,\dots,x^n\rangle/[\{ w_{i\>  i_1\dots i_N} x^{i_1}\otimes \dots \otimes x^{i_N}\}]
\end{equation}
where $w_{i_0  i_1\dots i_N}=w(e_{i_0},\dots e_{i_N})$ are the components of a 3-regular $(N+1)$-linear form $w$ on $\mathbb K^n$ with $N\geq 2$.
\end{proposition}

We write for the relations in Proposition \ref{D3}
\begin{equation}\label{dervol}
\partial_iw = w_{ii_1\dots i_N} x^{i_1}\otimes \dots \otimes x^{i_N}
\end{equation}
for $i\in\{1,\dots,n\}$. This is a sort of derivative of $w$ identified with
\begin{equation}\label{vol}
w=w_{i_0\dots i_N} x^{i_0}\otimes \dots \otimes x^{i_N}
\end{equation}
which is the generalization of the volume form and is referred to as the {\sl superpotential} \cite{vdb:2015}, \cite{boc:2008}, \cite{boc-sch-wem:2010}. For the interpretation in term of noncommutative volume see in \cite{mdv:2007} Proposition 10. \\
 
 Propositions \ref{D2} and \ref{D3} generalize to higher global dimensions but one has then to assume the $N$-Koszul property and one has to take higher order derivations of the superpotential (a preregular multilinear form) to write the relations, see in  \cite{mdv:2005} Theorem 4.3 and in \cite{mdv:2007} Theorem 11 and \cite{boc-sch-wem:2010} for the generalization to the quiver case.\\

\section{Exterior 3-forms}\label{Ef}

Let $(e_i) _{i\in\{1,\dots,n\}}$ be the canonical basis of $\mathbb K^n$ and let us equip $\mathbb K^n$ with the unique scalar product for which the canonical basis is orthonormal that is for which one has
\[
(e_i,e_j)=\delta_{ij}
\]
for $i,j\in\{1,\dots,n\}$. In the following, we assume that $n\geq3$ since we are interested in exterior 3-forms.\\

Let $\alpha$ be an exterior 3-form on $\mathbb K^n$. To $\alpha\in \wedge^3\mathbb K^n$ one associates $n$ endormorphisms $A_k$ of $\mathbb K^n$ by setting
\begin{equation}\label{mA}
A_ke_i=\sum_j\alpha_{jki}e_j=(A_k)^j_ie_j
\end{equation}
for $k,i\in \{1,\dots, n\}$, where $\alpha_{jki}=\alpha(e_j,e_k,e_i)$ are the components of $\alpha$. Since they are antisymmetric, the $A_k$ are $n$ elements of the Lie algebra $\fracso(n,\mathbb K)$ of $SO(n,\mathbb K)$. The 3-form
 $\alpha$  is cyclic and therefore $\alpha$  is preregular if and only if it is {\sl nondegenerate}, that is
\begin{equation}\label{ndeg}
i_X(\alpha)=0\Rightarrow X=0
\end{equation}
for $X\in \mathbb K^n$ where the 2-form $i_X(\alpha)$ is defined by
\begin{equation}\label{int}
i_X(\alpha)(Y,Z)=\alpha(X,Y,Z)
\end{equation}
for $Y,Z\in \mathbb K^n$. Condition (\ref{ndeg}) reads 
\begin{equation}\label{nd}
A_kX=0,\>\>\> \forall k\Rightarrow X=0
\end{equation}
in terms of the $A_k$. Finally, $\alpha$ is 3-regular if and only if 
\begin{equation}\label{3reg}
MA_k=A_kN,\>\>\> \forall k\Rightarrow M=N=\lambda\bbbone
\end{equation}
for $M,N\in \End(\mathbb K^n)$ with $\lambda\in \mathbb K$. Thus if $\alpha$ is 3-regular the system $(A_k)$ is irreductible (in view of the Schur lemma) which implies in particular that the $A_k\in \fracso(n,\mathbb K)$ generate the algebra $\End(\mathbb K^n)$ of the endomorphisms of $\mathbb K^n$ which contains as subspace the whole Lie algebra $\fracso(n,\mathbb K)$  of the antisymmetric endormorphisms of $\mathbb K^n$. However one should be aware of the fact that this does not mean that $\fracso(n,\mathbb K)$ is generated {\sl as Lie algebra} by the $A_k$. For instance the Lie algebra $\fracg_2(\mathbb K)$ has an irreducible representation in $\mathbb K^7$ although it is a proper Lie sub-algebra of $\fracso(7,\mathbb K)$.  In our general setting there is no rule as shown by the 2 families of examples given below.\\

In the next section we shall define exterior 3-forms $\alpha^{(p)}$ in dimensions $2p+1$ which are 3-regular and such that the corresponding $A^{(p)}_k\in \fracso(2p+1,\mathbb K)$ generate the whole Lie algebra $\fracso(2p+1,\mathbb K)$ of antisymmetric $(2p+1)\times (2p+1)$ matrices with coefficient in $\mathbb K$.\\

On the other hand, let $\alpha_{ijk}$ be the structure constants of the compact real form of a simple complex Lie algebra $\fracg$ of dimension $n$ in its standard orthonormal basis. Then the $\alpha_{ijk}$ define a real 3-regular 3-form $\alpha\in \wedge^3\mathbb R^n$ and the associated $A_k$ span the adjoint representation of $\fracg$ that is the corresponding Lie-sub-algebra $\pfraca$ of $\fracso(n,\mathbb R)$. The Jacobi identity reads then
\begin{equation}\label{mvAa}
(A_k)^j_{i_1} \alpha_{j i_2 i_3} + (A_k)^j_{i_2}\alpha_{i_1ji_3} + (A_k)^j_{i_3}  \alpha_{i_1 i_2j}=0
\end{equation}
for any $k\in \{1,\dots,n\}$, which means that the Lie sub-algebra $\pfraca$ of $\fracso(n,\mathbb R)$ generated by the $A_k$ preserves the exterior 3-form $\alpha$. This is clearly not the case for the whole Lie algebra $\fracso(n,\mathbb R)$ except for the dimension $n=3$, i.e. for $\fracg=\pfraca_1$ which corresponds to the exterior 3-form $\alpha^{(1)}$. The 3-regularity of $\alpha$ follows from the irreductibillity of the adjoint representation of $\fracg$.\\

\noindent {\bf Remark}\label{d4}. The dimension 4 is an exception since any 3-form $\alpha\in \wedge^3\mathbb K^4$ is of the form $\alpha=i_X(vol_4)$ for some $X\in \mathbb K^4$ where the volume form $vol_4\in \wedge^4\mathbb K^4$ is defined by $vol_4(e_1,e_2,e_3,e_4)=1$ which implies $vol_4(e_i,e_j,e_k, e_\ell)=\varepsilon_{ijk\ell}$. Theferore an element $\alpha=i_X(vol_4)\in \wedge^3\mathbb K^4$ is degenerate since then $i_X(\alpha)=0$ for $X\in \mathbb K^4$ with $X\not= 0$ if $\alpha\not=0$ and thus no non trivial $\alpha\in\wedge^3\mathbb K^4$ can be 3-regular. This has a counterpart on the side of the Lie algebra $\fracso(4,\mathbb K)$ which is not simple since $\fracso(4,\mathbb K)=\fracso(3,\mathbb K)\oplus \fracso(3,\mathbb K)$ while $\fracso(3,\mathbb K)$ and $\fracso(n,\mathbb K)$ for $n\geq 5$ are simple Lie algebras.\\

The first interesting non trivial case can only occur in dimension 5 and it turns out that it consists in the $GL(5,\mathbb K)$ orbit of the exterior 3-form $\alpha^{(2)}$.

\section{Exterior 3-forms as superpotentials}\label{Sf}

In this paper we are interested in the regularity of quadratic algebras with relations generated by exterior 3-forms that is by algebras of the form
\begin{equation}\label{A}
\cala=\mathbb K\langle x^1,\dots, x^n\rangle/[\{\alpha_{ijk} x^j\otimes x^k\}]
\end{equation}
where the $\alpha_{ijk}$ are the components of an exterior 3-form $\alpha$, i.e. are completely antisymmetric. Since by Proposition \ref{D3} we know that the 3-regularity of $\alpha$ is a necessary condition for the regularity of $\cala$, we are led to introduce a list of 3-regular 3-forms $\alpha$ in $\mathbb K^n$ and to study the regularity of the corresponding quadratic algebras given by (\ref{A}). We will use freely the following lemma.  

\begin{lemma}\label{kos-reg}
Let $\alpha$ be an exterior 3-regular 3-form on $\mathbb K^n$. Then the quadratic algebra
\[
\cala=\mathbb K\langle x^1,\dots, x^n\rangle/[\{\partial_i\alpha\}]
\]
is a domain and the following statements are equivalent : \\
$(\imath)$	$\cala$ is Koszul,\\
$(\imath\imath)$ $\cala$ is regular,\\
$(\imath \imath\imath)$ $\cala$ is 3-Calabi-Yau.
\end{lemma}
\noindent  {\bf Sketch of proof.} The algebra $\cala$ is a domain because it is the universal enveloping algebra of a Lie algebra since the relations involve only commutators \cite{jac:1962}. The equivalence $(\imath) \Leftrightarrow (\imath\imath)$  follows directly from Proposition 16  of 
\cite{mdv:2007}. Concerning the last equivalence, since $\alpha$ is cyclic $\cala$ is Calabi-Yau \cite{gin:2006} whenever it is regular \cite{boc:2008} . $\square$\\

Let $\alpha\in \wedge^3\mathbb K^n$ and $\cala$ be the quadratic algebra defined by (\ref{A}) . One defines an antisymmetric $n\times n$-matrix  $A(x)$ with entries in $\cala$ by setting
\begin{equation}\label{Ax}
A(x)=A_k x^k
\end{equation}
where the $A_k$ are the matrices given by (\ref{mA}), that is by $(A_k)^i_j=\alpha_{ikj}$.
The relations $\partial_i\alpha$ of $\cala$ can be expressed as
\begin{equation}\label{relA}
\left(
\begin{array}{c}
\partial_1 \alpha\\
\vdots\\
\partial_n\alpha
\end{array}
\right)
= A(x)
\left(
\begin{array}{c}
x^1\\
\vdots\\
x^n
\end{array}
\right)
\end{equation}
or simply by $\partial x=A(x)x$ where $\partial x$ and $x$ denote the corresponding columns. Let us assume now that $\alpha$ is 3-regular, then the augmented Koszul complex of $\cala$ reads as follows from \cite{mdv:2007}
\begin{equation}\label{KA}
0\longrightarrow \cala \stackrel{x^t}{\longrightarrow}\cala^n \stackrel{A(x)}{\longrightarrow}\cala^n\stackrel{x}{\longrightarrow} \cala\stackrel{\varepsilon}{\longrightarrow} \mathbb K \longrightarrow 0
\end{equation}
where $x^t$ means multiplication in $\cala$ by the transposed $x^t$ of the column $x$, $A(x)$ means multiplication in $\cala$ of the line $\cala^n$ of elements of $\cala$ with the $n\times n$-matrix $A(x)$ of elements of $\cala$,  $x$ means multiplication in $\cala$ of the line $\cala^n$ of elements of $\cala$ by the column $x$ while $\varepsilon$ is the projection onto degree 0.\\
The exactness of the sequence 
\[
0\longrightarrow \cala \stackrel{x^t}{\longrightarrow} \cala^n
\]
follows from the fact that $\cala$ is a domain while the exactness of the sequence
\[
\cala^n\stackrel{A(x)}{\longrightarrow} \cala^n \stackrel{x}{\longrightarrow} \cala\stackrel{\varepsilon}{\longrightarrow}\mathbb K \longrightarrow 0
\]
just expresses the definition of $\cala$ by generators and relations. Therefore $\cala$ is Koszul if and only if the sequence
\begin{equation}\label{exRel}
\cala\stackrel{x^t}{\longrightarrow} \cala^n \stackrel{A(x)}{\longrightarrow} \cala^n
\end{equation}
is exact.\\

It is clear that the quadratic algebras $\cala$ of the form (\ref{A}) do only depend up to isomorphism on the 3-forms $\alpha$ up to the $GL(n,\mathbb K)$ action. Thus we are only interested in the $GL(n,\mathbb K)$ orbits of 3-forms in $\mathbb K^n$. Notice that the 3-regularity of a 3-form is an invariant notion as well as the regularity of the corresponding algebra. In Section \ref{3F} we shall analyse the orbits of nondegenerate exterior 3-forms for $\mathbb K=\mathbb C$ by choosing a convenient element in each orbit for dimension $n$ with $3\leq n\leq 7$, using the results of 
\cite{coh-hel:1988}. We will then select the 3-regular orbits and corresponding 3-regular representative exterior 3-forms. The regularity of the corresponding quadratic algebras will be investigated in Section \ref{RF}.\\

We now analyse a particular family of exterior 3-forms $\alpha^{(p)}\in \wedge^3\mathbb K^{2p+1}$ and show that the corresponding quadratic algebras $\cala^{(p)}$ are regular of global dimension 3. Let us define $\alpha^{(p)}$ by 
\begin{equation}\label{alf}
\alpha^{(p)}=\sum^p_{m=1}\theta^m \wedge \theta^{m+p}\wedge \theta^{2p+1}=\alpha_{ijk}\theta^i\otimes \theta^j\otimes \theta^k
\end{equation}
where $(\theta^i)$ is the dual basis of the basis $(e_j)$ of $\mathbb K^{2p+1}$. It is easily verified that $\alpha^{(p)}$ is 3-regular. The corresponding superpotential is $\alpha^{(p)}=\alpha_{ijk}x^i\otimes x^j\otimes x^k$ and the relations $\partial_i\alpha^{(p)}$ of the corresponding quadratic algebra $\cala^{(p)}$ read
\begin{equation}\label{rp1}
\begin{array}{l}
\partial_m\alpha^{(p)}=x^{m+p}\otimes x^{2p+1}-x^{2p+1}\otimes x^{m+p}\\
\partial_{m+p} \alpha^{(p)}=x^{2p+1}\otimes x^m - x^m\otimes x^{2p+1}
\end{array}
\end{equation}
for $1\leq m\leq p$ and
\begin{equation}\label{rp2}
\partial_{2p+1} \alpha^{(p)}=\sum^p_{r=1} x^r\otimes x^{r+p} - x^{r+p}\otimes x^r
\end{equation}
and imply that $x^{2p+1}$ is in the center of $\cala^{(p)}$ in view of (\ref{rp1}) and that the $x^\ell$ for $\ell\in \{1,\dots,2p\}$ span the quadratic sub-algebra
\begin{equation}\label{Bp}
\calb^{(p)}= \mathbb K \langle x^1,\dots, x^{2p}\rangle / [\partial_{2p+1}\alpha^{(p)}]
\end{equation}
of $\cala^{(p)}$. Since $\partial_{2p+1}\alpha^{(p)}$ given by (\ref{rp2}) is obviously the superpotential corresponding to a nondegenerate bilinear form on $\mathbb K^{2p}$, $\calb^{(p)}$ is a regular quadratic algebra of global dimension 2. Thus one has
\begin{equation}\label{Ap}
\cala^{(p)}= \calb^{(p)} \otimes  \mathbb K[x^{2p+1}]
\end{equation}
which implies that $\cala^{(p)}$ is regular of global dimension 3.\\

Let us now show that the $(2p+1)\times (2p+1)$ antisymmetric matrices $A^{(p)}_k\in \fracso(2p+1,\mathbb K)$ associated (via (\ref{mA})) to $\alpha^{(p)}\in \wedge^3\mathbb K^{2p+1}$ generate as Lie algebra the whole Lie algebra $\fracso(2p+1,\mathbb K)$. One defines
\begin{equation}\label{Au}
A(u)=\left (
\begin{array}{ccc}
0_p & -u^{2p+1}\bbbone_p & \begin{array}{c}u^{2p} \\ \vdots\\ u^{p+1} \end{array}\\
\\
u^{2p+1}\bbbone_p & 0_p & \begin{array}{c} -u^1\\ \vdots \\  -u^p \end{array}\\
\\
-u^{2p}\dots -u^{p+1} & u^1\dots u^p & 0
\end{array}
\right)
\end{equation}
by setting $A(u)= A^{(p)}_k u^k$ for $u\in \mathbb K^{2p+1}$. Then by computing the commutator $[A(u),A(v)]$ for $u,v\in \mathbb K^{2p+1}$, one verifies that the components of the antisymmetric matrice $[A(u),A(v)]$ are just a permutation of the components of the exterior product $u\wedge v$. It is clear that the linear span of the $u\wedge v$ for $u,v\in \mathbb K^{2p+1}$ is the whole Lie algebra $\fracso(2p+1,\mathbb K)$ of antisymmetric $(2p+1)\times (2p+1)$-matrices, which is therefore also the case for the $[A(u),A(v)]$ for $u,v\in \mathbb K^{2p+1}$.\\

We summarize the above results concerning the $\alpha^{(p)}$ by the following proposition.

\begin{proposition}\label{Aap}
Let $\alpha^{(p)}\in \wedge^3\mathbb K^{2p+1}$ be defined by (\ref{alf}) and let $\cala^{(p)}$ be the associated quadratic algebra as in (\ref{A}) and $A^{(p)}_k$ be the corresponding endomorphisms of $\mathbb K^{2p+1}$ as in (\ref{mA}). Then one has the following .
\begin{enumerate}
\item The quadratic algebra $\cala^{(p)}$ is regular of global dimension 3.
\item The $A^{(p)}_k$, $k\in \{1,\dots, 2p+1\}$ generate the Lie-algebra $\fracso(2p+1,\mathbb K)$.
\end{enumerate}
\end{proposition}

\section{Exterior 3-forms in $\mathbb K^n$ with $3\leq n\leq 7$}\label{3F}

In this section $\mathbb K$ is algebraically closed. We shall use the results of \cite{coh-hel:1988} for $3\leq n\leq 7$ 
and the choices of representative elements in $\wedge^3\mathbb K^n$ for the $GL(n,\mathbb K)$-orbits of nondegenerate (i.e. preregular) exterior 3-forms, see also in  \cite{gur:1964}. However the analysis of the regularity of these elements is independent of the assumption that $\mathbb K$ is algebraically closed.\\

{\bf In dimension $n=3$}, the only nontrivial orbit is the one of $\alpha^{(1)}=\theta^1\wedge \theta^2\wedge\theta^3$. The corresponding superpotential is
\begin{equation}\label{sp1}
\alpha^{(1)}=\varepsilon_{ijk} x^i\otimes x^j\otimes x^k
\end{equation}
and the quadratic algebra $\cala^{(1)}$ is the commutative algebra
\begin{equation}\label{cala1}
\cala^{(1)}=\mathbb K[x^1,x^2,x^3]
\end{equation}
of polynomials in the $x^i$, $i\in \{1,2,3\}$. By setting $A^{(1)}(x^1,x^2,x^3)=A^{(1)}_kx^k$ one has
\begin{equation}\label{mat1}
A^{(1)}(x^1,x^2,x^3)=\left(\begin{array}{ccc}
0 & -x^3 & x^2\\
x^3 & 0 & -x^1\\
-x^2 & x^1 & 0
\end{array}
\right)
\end{equation}
for the corresponding matrix.\\

{\bf In dimension $n=4$}, all exterior 3-forms are degenerate.\\

{\bf In dimension $n=5$}, the only orbit of nondegenerate exterior 3-forms is the one of $\alpha^{(2)}=(\theta^1\wedge \theta^3 +\theta^2\wedge \theta^4)\wedge \theta^5=\alpha^{(2)}_{ijk}\theta^i\wedge \theta^j\wedge \theta^k$ which is 3-regular. Moreover the associated quadratic algebra $\cala^{(2)}$ is regular of global dimension 3 (see in the last section).\\

{\bf In dimension $n=6$}, there are 2 orbits of nondegenerate exterior 3-forms : namely the orbit of $\gamma=\theta^1\wedge \theta^2\wedge \theta^3+\theta^4\wedge \theta^5\wedge \theta^6$ and the orbit of $\omega=\theta^1\wedge \theta^2 \wedge \theta^6+\theta^3 \wedge \theta^1 \wedge \theta^5 + \theta^2 \wedge \theta^3 \wedge \theta^4$.
None of these 3-forms is 3-regular so the associated quadratic algebras $\cala_\gamma$ and $\cala_\omega$ cannot be regular. Indeed $\gamma$ is disconnected that is the sum of two exterior 3-forms in two different spaces $\mathbb K^3$ with basis ($e_1,e_2,e_3$) and ($e_4,e_5,e_6$) respectively, so the commutant of the corresponding  
\[
A_\gamma (x) =\left(\begin{array}{cc}
A^{(1)}(x^1,x^2,x^3)& 0_3\\
0_3 & A^{(1)} (x^4,x^5, X^6) 
\end{array}
\right)
\]
contains the matrices 
\[
\left(\begin{array}{cc}
\lambda\bbbone_3 & 0_3\\
0_3 & \mu\bbbone_3
\end{array}
\right)
\]
for any $\lambda, \mu\in \mathbb K$ which implies that $\gamma$ is not 3-regular, while for $\omega$, one has
\[
A_\omega(x)=\left(
\begin{array}{cc}
A^{(1)}(x^4,x^5,x^6)  & -A^{(1)}(x^1,x^2, x^3)^t\\
A^{(1)}(x^1,x^2,x^3) & 0_3
\end{array}
\right)
\]
where $A^{(1)}$ is defined by (\ref{mat1}),  so one has for any $\mu\in \mathbb K$
\[
\left(
\begin{array}{cc}
0_3 & \mu\bbbone_3\\
0_3 & 0_3
\end{array}
\right)
A_\omega = A_\omega
\left(
\begin{array}{cc}
0_3 & 0_3\\
\mu\bbbone_3 & 0_3
\end{array}
\right)
\]
which implies that $\omega$ is not 3-regular.\\

{\bf In dimension $n=7$} , the situation is much more elaborated. There are five orbits of nondegenerate exterior 3-forms namely the orbits of the following exterior 3-forms:
\[
\begin{array}{l}
\rho=\theta^1 \wedge \theta^2 \wedge \theta^3 + \theta^2\wedge \theta^4\wedge \theta^6 +\theta^3\wedge \theta^5 \wedge \theta^7\\
\beta=(\theta^1\wedge \theta^3+\theta^2\wedge \theta^4)\wedge \theta^5 +  \theta^1\wedge \theta^7 \wedge \theta^2 +\theta^3\wedge \theta^6 \wedge \theta^4\\
\alpha^{(3)}=(\theta^1\wedge\theta^4+\theta^2\wedge \theta^5 + \theta^3\wedge\theta^6)\wedge \theta^7\\
\alpha^{(3)\prime}=\alpha^{(3)}+\theta^1\wedge \theta^2 \wedge \theta^3\\
\alpha^{(3)\prime\prime}=\alpha^{(3)}+\theta^1\wedge \theta^2 \wedge \theta^3 +\theta^4\wedge \theta^5 \wedge \theta^6
\end{array}
\]
where we have used the notation of Section 3 for $\alpha^{(3)}$. It is worth noticing here that one can also express $\beta$ as
\[
\beta =\alpha^{(2)}+\theta^1\wedge \theta^7\wedge \theta^2+\theta^3\wedge \theta^6\wedge \theta^4
\]
with the notation of Section 3. The correspondence with the notations of \cite{coh-hel:1988} is $\alpha^{(1)}=f_1, \alpha^{(2)}=f_2, \gamma=f_3, \omega=f_4,\rho=f_5, \beta=f_7, \alpha^{(3)}=f_8,\alpha^{(3)\prime}
=f_6$ and $\alpha^{(3)\prime\prime}=f_9$.\\

The exterior 3-forms $\rho$ and $\beta$ are 3-regular while $\alpha^{(3)},\alpha^{(3)\prime}$ and $\alpha^{(3)\prime\prime}$ span an affine plane
\begin{equation}
\alpha^{(t_0,t_1,t_2)}=t_0 \alpha^{(3)}+t_1\alpha^{(3)\prime}+t_2\alpha^{(3)\prime\prime}, \>\> t_0+t_1+t_2=1
\end{equation}
of 3-regular exterior 3-forms to which are associated a corresponding family $\cala^{(t_0,t_1,t_2)}$ of quadratic 3-Calabi-Yau algebras of the type investigated in \cite{he-fvo-zha:2015}. Namely the $\cala^{(t_0,t_1,t_2)}$ are the cross-products of the 2-Calabi-Yau algebra $\calb^{(3)}$ defined by (\ref{Bp}) with the derivations $t_1\delta_1+t_2\delta_2$, (i.e. are the Ore extensions of $\calb^{(3)}$  associated with data $(I,t_1\delta_1+t_2\delta_2)$). Notice that $\alpha^{(3)\prime\prime}$ is in the orbit of the generic 3-forms so the associated quadratic algebra is isomorphic to the algebra	 with $G_2$-symmetry defined in \cite{smi:2011}.\\

For $\alpha=\rho$ the matrix $A(x)$ given by  (\ref{Ax}) reads

\begin{equation}\label{Arho}
A_\rho(x)=\left(
\begin{array}{ccccccc}
0 & -x^3 & x^2 & 0 & 0 & 0 & 0\\
x^3 & 0 & -x^1 & -x^6 & 0 & x^4 & 0\\
-x^2 & x^1 & 0 & 0 & -x^7 & 0 & x^5\\
0 & x^6 & 0 & 0 & 0 & -x^2 & 0\\
0 & 0 & x^7 & 0 & 0  & 0 & -x^3\\
0 & -x^4 & 0 & x^2 & 0 & 0 & 0\\
0 & 0 & -x^5 & 0  & x^3 & 0 & 0
\end{array}
\right)
\end{equation}
and the relations of the associated quadratic algebra $\cala_\rho$ read
\[
\begin{array}{lll}
\partial_1\rho & = & x^2\otimes x^3-x^3\otimes x^2\\
\partial_2\rho & = & x^3\otimes x^1-x^1\otimes x^3  + x^4\otimes x^6 - x^6\otimes x^4\\
\partial_3 \rho & = & x^1\otimes x^2 - x^2\otimes x^1 + x^5\otimes x^7-x^7\otimes x^5\\
\partial_4\rho & = & x^6\otimes x^2 - x^2\otimes x^6\\
\partial_5\rho & = & x^7\otimes x^3 - x^3\otimes x^7\\
\partial_6\rho & = & x^2\otimes x^4 - x^4\otimes x^2\\
\partial_7\rho & = & x^3\otimes x^5 - x^5 \otimes x^3
\end{array}
\]
One can show by direct computation (or by using computer) that for matrices $M,N\in M_7(\mathbb K)$, the relation
\[
MA_\rho(x) = A_\rho(x) N
\]
implies $M=N=\lambda\bbbone_7$, for some $\lambda\in \mathbb K$ so the exterior 3-form $\rho$ is 3-regular.\\

For $\alpha=\beta$, one has for the corresponding matrix $A(x)=A_\beta(x)$

\begin{equation}\label{AB7}
A_\beta(x) = \left(
\begin{array}{ccccccc}
0&x^7&-x^5&0&x^3&0&-x^2\\
-x^7&0&0&-x^5&x^4&0&x^1\\
x^5&0&0&x^6&-x^1&-x^4&0\\
0&x^5&-x^6&0&-x^2&x^3&0\\
-x^3&-x^4&x^1&x^2&0&0&0\\
0&0&x^4&-x^3&0 &0&0\\
x^2&-x^1&0&0&0&0&0
\end{array}
\right)
\end{equation}
and the relations of the associated quadratic algebra $\cala_\beta$ read
\[
\begin{array}{l}
\partial_1\beta = x^3 \otimes x^5-x^5\otimes x^3 + x^7\otimes x^2 - x^2\otimes x^7\\
\partial_2\beta = x^4 \otimes x^5-x^5\otimes x^4 + x^1\otimes x^7 - x^7\otimes x^1\\
\partial_3\beta = x^5 \otimes x^1-x^1\otimes x^5 + x^6\otimes x^4 - x^4\otimes x^6\\
\partial_4\beta = x^5 \otimes x^2-x^2\otimes x^5 + x^3\otimes x^6 - x^6\otimes x^3\\
\partial_5\beta = x^1 \otimes x^3-x^3\otimes x^1 + x^2\otimes x^4 - x^4\otimes x^2\\
\partial_6\beta=x^4\otimes x^3 - x^3\otimes x^4\\
\partial_7\beta=x^2\otimes x^1 - x^1\otimes x^2
\end{array}
\]

Thus it follows that $\cala_\beta$ has 2 non intersecting sub-algebras namely the quadratic sub-algebra $\calb$ generated by $x^1,x^2, x^3, x^4$ with relations $\partial_5\beta, \partial_6\beta, \partial_7\beta$ and the sub-algebra $\calc=\mathbb K\langle x^5,x^6, x^7\rangle$ freely generated by $x^5,x^6,x^7$ which is isomorphic to the tensor algebra $T(\mathbb K^3)$ and is a Koszul algebra of global dimension one. Concerning the sub-algebra $\calb$ it is a Koszul algebra of global dimension 2. Indeed by introducing the $3\times4$-matrix $B(x)$ with coefficients in $\calb$ defined by
\begin{equation}\label{B7}
B(x) =\left(
\begin{array}{cccc}
-x^3&-x^4&x^1&x^2\\
0&0&x^4&-x^3\\
x^2&-x^1&0&0
\end{array}
\right)
\end{equation}
the augmented Koszul complex of $\calb$ reads
\[
0\longrightarrow \calb^3 \stackrel{B(x)}{\longrightarrow} \calb^4\stackrel{x}{\longrightarrow}\calb\stackrel{\varepsilon}{ \longrightarrow}\mathbb K\longrightarrow 0
\]
and, by using the fact that $\calb$ is a domain since it is a universal enveloping algebra of a Lie algebra, it is easy to show that this is an exact sequence, that is that the mapping $\calb^3\rightarrow \calb^4$ is injective. Thus $\cala_\beta$ is generated by the two non intersecting sub-algebra $\calb$ and $\calc$ which are Koszul of global dimensions $\gldim(\calb)=2$ and $\gldim(\calc)=1$. The relations $\partial_1\beta, \partial_2\beta, \partial_3\beta, \partial_4\beta$ give commutation relations between the generators of $\calb$ and generators of $\calc=\mathbb K\langle x^5,x^6,x^7\rangle$.\\

Finally one can show by using computer that for matrices $M,N\in M_7(\mathbb K)$ the relation
\[
MA_\beta(x)=A_\beta(x)N
\]
implies $M=N=\lambda\bbbone_7$ for some $\lambda\in\mathbb K$ which means that $\beta$ is 3-regular.\\

We now discuss the cases of $\alpha^{(3)\prime}$ and $\alpha^{(3)\prime\prime}$.\\

For $\alpha^{(3)}$ this is the discussion of Section \ref{Sf} for $p=3$. One has $\cala^{(3)}=\calb^{(3)}\otimes \mathbb K[x^7]$ which is the cross-product of $\calb^{(3)}$ with the trivial derivation $\delta_0=0$.\\

For $\alpha^{(3)\prime}$, one has
\[
A^{(3)\prime}(x)=\left (
\begin{array}{ccccccc}
0&-x^3&x^2&-x^7&0&0&x^4\\
x^3&0&-x^1&0&-x^7&0&x^5\\
-x^2&x^1&0&0&0&-x^7&x^6\\
x^7&0&0&0&0&0&-x^1\\
0&x^7&0&0&0&0&-x^2\\
0&0&x^7&0&0&0&-x^3\\
-x^4 & -x^5& -x^6 &x^1 & x^2 & x^3&0
\end{array}
\right)
\]
for the corresponding matrix $A(x)=A^{(3)\prime}(x)$ and the relations of the associated quadratic algebra $\cala^{(3)\prime}$ read
\[
\begin{array}{l}
\partial_1\alpha^{(3)\prime}  =  x^4\otimes x^7 -x^7\otimes x^4 + x^2\otimes x^3 -x^3\otimes x^2\\
\partial_2\alpha^{(3)\prime}  =  x^5\otimes x^7 -x^7\otimes x^5 + x^3\otimes x^1 -x^1\otimes x^3\\
\partial_3\alpha^{(3)\prime}  =  x^6\otimes x^7 -x^7\otimes x^6 + x^1\otimes x^2 -x^2\otimes x^1\\
\partial_4\alpha^{(3)\prime} = \partial_4\alpha^{(3)}  =  x^7\otimes x^1-x^1\otimes x^7\\
\partial_5\alpha^{(3)\prime} = \partial_5\alpha^{(3)}  =  x^7\otimes x^2-x^2\otimes x^7\\
\partial_6\alpha^{(3)\prime}= \partial_6\alpha^{(3)}  =  x^7\otimes x^3-x^3\otimes x^7\\
\partial_7\alpha^{(3)\prime}=x^1\otimes x^4-x^4\otimes x^1+x^2\otimes x^5-x^5\otimes x^2+x^3\otimes x^6-x^6\otimes x^3\\
=\partial_7\alpha^{(3)}
\end{array}
\]
Now let $\calb^{(3)}$ be the quadratic algebra $\calb^{(p)}$ of Section \ref{Sf} for $p=3$, that is
\[
\calb^{(3)}=\mathbb K\langle x^1,\dots,x^6\rangle /[\partial_7 \alpha^{(3)}]
\]
which is a regular algebra of global dimension 2, then
\[
\delta_1x^1=x^2\otimes x^3-x^3\otimes x^2,\> \delta_1x^2=x^3\otimes x^1-x^1\otimes x^3,\> \delta_1x^3=x^1\otimes x^2-x^2\otimes x^1,
\]
$\delta_1x^4=0$, $\delta_1x^5=0$ and $\delta_1 x^6=0$ defines a derivation of degree 1 of $\mathbb K\langle x^1,\dots,x^6\rangle$ which satisfies
\[
\delta_1(\partial_7\alpha^{(3)})=0
\]
in view of the associativity of the tensor product. It follows that $\delta_1$ defines a derivation of degree 1, again denoted by $\delta_1$, of $\calb^{(3)}$. Then the relations of $\cala^{(3)\prime}$ imply that $\cala^{(3)\prime}$ is the cross-product 
\[
\cala^{(3)\prime}=\calb^{(3)}\rtimes_{\delta_1}\mathbb K[x^7]
\]
of $\calb^{3)}$ with the derivation $\delta_1$, the corresponding new generator being $x^7$, $(\delta_1=ad(x^7))$.\\

One proceeds similarily for $\alpha^{(3)\prime\prime}$, one has
\[
A^{(3)\prime\prime}(x)=\left (
\begin{array}{ccccccc}
0&-x^3&x^2&-x^7&0&0&x^4\\
x^3&0&-x^1&0&-x^7&0&x^5\\
-x^2&x^1&0&0&0&-x^7&x^6\\
x^7&0&0&0&-x^6&x^5&-x^1\\
0&x^7&0&x^6&0&-x^4&-x^2\\
0&0&x^7&-x^5&x^4&0&-x^3\\
-x^4 & -x^5& -x^6 &x^1 & x^2 & x^3&0
\end{array}
\right)
\]
while the relations of $\cala^{(3)\prime\prime}$ read 
\[
\begin{array}{l}
\partial_1\alpha^{(3)\prime\prime}= \partial_1\alpha^{(3)\prime}  =  x^4\otimes x^7 -x^7\otimes x^4 + x^2\otimes x^3 -x^3\otimes x^2\\
\partial_2\alpha^{(3)\prime\prime}= \partial_2\alpha^{(3)\prime}  =  x^5\otimes x^7 -x^7\otimes x^5 + x^3\otimes x^1 -x^1\otimes x^3\\
\partial_3\alpha^{(3)\prime\prime}= \partial_3\alpha^{(3)\prime}  =  x^6\otimes x^7 -x^7\otimes x^6 + x^1\otimes x^2 -x^2\otimes x^1\\
\partial_4\alpha^{(3)\prime\prime} =  x^7\otimes x^1-x^1\otimes x^7 +x^5\otimes x^6-x^6\otimes x^5\\
\partial_5\alpha^{(3)\prime\prime}  =  x^7\otimes x^2-x^2\otimes x^7 + x^6\otimes x^4-x^4\otimes x^6\\
\partial_6\alpha^{(3)\prime\prime} = x^7\otimes x^3-x^3\otimes x^7+ x^4\otimes x^5-x^5\otimes x^4\\
\partial_7\alpha^{(3)\prime\prime}=\partial_7\alpha^{(3)\prime}=\partial_7\alpha^{(3)}
\end{array}
\]
One defines a derivation $\delta_2$ of degree 1 of $\mathbb K\langle x^1,\dots, x^6\rangle$ by setting\\

\noindent $\delta_2x^1=x^2\otimes x^3-x^3\otimes x^2,\> \delta_2x^2=x^3\otimes x^1-x^1\otimes x^3,\> \delta_2x^3=x^1\otimes x^2-x^2\otimes x^1,\\ 
\delta_2x^4=x^5\otimes x^6-x^6\otimes x^5,\delta_2x^5=x^6\otimes x^4-x^4\otimes x^6,\> \delta_2x^6=x^4\otimes x^5-x^5\otimes x^2$\\
 
which satisfies again
\[
\delta_2(\partial_7\alpha^{(3)})=0
\]
in view of the associativity of the tensor product. Therefore $\delta_2$ passes to the quotient and defines a derivation of degree 1, again denoted by $\delta_2$, of the regular algebra $\calb^{(3)}$. Thus $\cala^{(3)\prime\prime}$ is the cross-product
\[
\cala^{(3)\prime\prime}=\calb^{(3)}\rtimes_{\delta_2}\mathbb K[x^7]
\]
of $\calb^{(3)}$ with the derivation $\delta_2$, the corresponding ``new" generator of degree 1 being $x^7$.\\

\noindent {\bf Remark}. Since the relations of all the quadratic algebras of this paper have relations defined in terms of commutators it follows that they are the universal enveloping algebras of graded Lie algebras generated in degree 1 as already pointed out. The cross-products with derivations considered above correspond to the universal enveloping algebra versions of semi-direct products of the graded Lie algebras with derivations of degree 1. Thus $\calb^{(3)}$ is the universal enveloping algebra of the corresponding graded Lie algebra $\pfracb^{(3)}$ and the $\delta_i, i \in \{0,1,2\}$ are in fact derivations of degree 1 of $\pfracb^{(3)}$ while $\cala^{(3)}, \cala^{(3)\prime}, \cala^{(3)\prime\prime}$ are the universal enveloping algebras of the semi-direct products of the Lie algebra $\pfracb^{(3)}$ with the derivations $\delta_0,\delta_1,\delta_2$.

\section{Regularity}\label{RF}

The cross-product with a derivation of degree 1 is a particular case of the graded Ore extension which preserves the Koszul property and even more, it preserves the $\calk_2$ property \cite{cas-she:2008} in a very strong sense \cite{pha:2012}. It follows that the quadratic algebras $\cala^{(3)}, \cala^{(3)\prime}, \cala^{(3)\prime\prime}$ associated with the exterior 3-forms $\alpha^{(3)}, \alpha^{(3)\prime}, \alpha^{(3)\prime\prime}\in \wedge^3\mathbb K^7$ are Koszul and therefore are regular in view of Lemma \ref{kos-reg}. This is also true for the quadratic algebras $\cala^{(p)}$ associated with the exterior 3-forms $\alpha^{(p)}\in \wedge^3\mathbb K^{2p+1}$  $(p\geq 1)$ as explained in the second part of Section \ref{Sf}. All these quadratic algebras are particular cases of the ones of \cite{he-fvo-zha:2015}.\\

It remains to discuss the case of the quadratic algebras $\cala_\rho$ and $\cala_\beta$ associated with the 3-forms $\rho$ and $\beta\in\wedge^3 \mathbb K^7$.\\

The matrix $A_\rho(x)$ is given by (\ref{Arho}) so the relations of the algebra $\cala_\rho$ read

\begin{equation}\label{rr1}
[x^2,x^3]=0
\end{equation}
\begin{equation}\label{rr2}
[x^3, x^1] + [x^4, x^6]  = 0 
\end{equation}
\begin{equation}\label{rr3}
[x^1,x^2] + [x^5,x^7]=0 
\end{equation}
\begin{equation}\label{rr4}
[x^6,x^2]=0
\end{equation}
\begin{equation}\label{rr5}
[x^7,x^3]=0
\end{equation}
\begin{equation}\label{rr6}
[x^2,x^4]=0
\end{equation}
\begin{equation}\label{rr7}
[x^3,x^5]=0\
\end{equation}
for the generators $x^k$ ($k\in \{1,\dots, 7\}$) of $\cala_\rho$.

\begin{proposition}\label{Regr}
The sequence 
\[
\cala_\rho \stackrel{x^t}{\longrightarrow} \cala^7_\rho \stackrel{A_\rho(x)}{\longrightarrow} \cala^7_\rho
\]
is exact.
\end{proposition}

\noindent \underbar{Proof}. It is sufficient to prove that
\begin{equation}\label{Kerr}
(a^1,a^2, a^3, a^4, a^5, a^6, a^7) A_\rho(x)=0
\end{equation}
for $(a^k)\in \cala^7_\rho$ implies that
\begin{equation}\label{axr}
a^k=ax^k\>\>\> \forall k\in \{1,\dots, 7\}
\end{equation}
for some $a\in \cala_\rho$.\\

Equation (\ref{Kerr}) reads by using (\ref{Arho})

\begin{equation}\label{ar1}
a^2x^3 - a^3 x^2=0
\end{equation}
\begin{equation}\label{ar2}
-a^1x^3 + a^3x^1 + a^4x^6-a^6x^4=0
\end{equation}
\begin{equation}\label{ar3}
a^1x^2 - a^2 x^1 + a^5 x^7-a^7x^5=0
\end{equation}
\begin{equation}\label{ar4}
-a^2x^6+a^6x^2=0
\end{equation}
\begin{equation}\label{ar5}
-a^3x^7+a^7x^3=0
\end{equation}
\begin{equation}\label{ar6}
a^2x^4-a^4x^2=0
\end{equation}
\begin{equation}\label{ar7}
a^3x^5-a^5x^3=0
\end{equation}

For the $a^k$ of degree 0 that is $a^k\subset \mathbb K$ ($\forall k\in \{1,\dots, 7\}$) it is clear that (\ref{Kerr})$\Rightarrow$ (\ref{axr}) with $a=0$ since the $x^k$ are linearly independent.\\

In order to prove the above proposition we first prove the following lemmas
\begin{lemma}\label{a23}
Assume that the $a^k\in \cala_\rho$ with $2\leq k\leq 7$ satisfy the relations (\ref{ar4}), (\ref{ar5}), (\ref{ar6}), (\ref{ar7}). Then $a^2=ax^2$ and $a^3=ax^3$ for some $a\in \cala_\rho$ imply $a^k=ax^k$ for any $k$ with $2\leq k\leq 7$.
\end{lemma}

\noindent \underbar{Proof}. Assume $a^2=ax^2$  and $a^3=ax^3$ then\\

(\ref{ar4}) reads $(a^6-ax^6)x^2=0$ in view of (\ref{rr4})\\

(\ref{ar5}) reads $(a^7-ax^7)x^3=0$ in view of (\ref{rr5})\\

(\ref{ar6}) reads $(a^4-ax^4)x^2=0$ in view of (\ref{rr6})\\

(\ref{ar7}) reads $(a^5-ax^5)x^3=0$ in view of (\ref{rr7}).\\

Since $\cala_\rho$ is a domain this implies $a^k=ax^k$ for $2\leq k\leq 7$. $\square$

\begin{lemma}\label{a123}
Assume now that the $a^k\subset \cala_\rho$ with $2\leq k\leq 7$ satisfy (\ref{ar1}), (\ref{ar4}), (\ref{ar5}), (\ref{ar6}),(\ref{ar7}). Then there is some $a\in \cala_\rho$ such that $a^k=ax^k$ for any  $2\leq k\leq 7$.
\end{lemma}

\noindent \underbar{Proof}. We proceed by induction on the degree $d$ of the $a^k$. This is clearly true for $d=0$, with $a=0$ in view of the linear independence of the $x^k$. Let us assume that this is true for $d\leq n$ and let the degree of the $a^k$ be $n+1$. The relation (\ref{ar1}) reads
\[
(0, -a^3, a^2,0,0,0,0)x= 0
\]
which implies that  there are $b^k\in \cala_p$ for $k\in \{1,\dots, 7\}$ such that
\begin{equation}\label{ba}
(b^1,\dots, b^7)A_\rho(x)=(0,-a^3, a^2, 0,0,0,0)
\end{equation}
since the exactness of the sequence 
\[
\cala^7_\rho \stackrel{A_\rho(x)}{\rightarrow} \cala^7_\rho \stackrel{x}{\rightarrow} \cala
\]
is equivalent to the presentation of $\cala_\rho$.

The relation (\ref{ba}) reads
\begin{equation}\label{ba23}
b^2 x^3 - b^3 x^2=0
\end{equation}
\begin{equation}\label{ba13}
- b^1 x^3 + b^3 x^1 + b^4 x^6 - b^6 x^4 = -a^3
\end{equation}
\begin{equation}\label{ba12}
b^1 x^2 - b^2 x^1 + b^5 x^7 - b^7 x^5 = a^2
\end{equation}
\begin{equation}\label{ba26}
- b^2 x^6 + b^6 x^2=0
\end{equation}
\begin{equation}\label{ba37}
-b^3 x^7 + b^7 x^3 =0
\end{equation}
\begin{equation}\label{ba24}
b^2x^4 - b^4 x^2=0
\end{equation}
\begin{equation}\label{ba35}
b^3x^5 - b^5 x^3=0.
\end{equation}
Now the $b^k$ are of degree $n$ and therefore for $k\geq 2$ one has  $b^k=bx^k$ for some $b\in \cala_\rho$ in view of Lemma \ref{a23} and of the induction assumption. Therefore, using (\ref{rr3}) and (\ref{rr2}),  relations (\ref{ba12}) and (\ref{ba13}) read $a^2=(b^1-bx^1)x^2$ and $a^3=(b^1-bx^1)x^3$ which in view of Lemma \ref{a23} implies
\[
a^k=ax^k\>\>\> \text{for}\>\> 2\leq k\leq 7
\]
with $a=b^1-bx^1$. $\square$\\

\noindent \underbar{End of proof of Proposition \ref{Regr}}.\\
Let $(a^1,\dots, a^7)\in \ker (A_\rho (x))$, we have proven that  $a^k=ax^k$ for $k\geq 2$ but then relations (\ref{ar3}) and (\ref{ar2}) read $(a^1-ax^1)x^2=0$ and $(a^1-ax^1)x^3=0$. Each one of the last relations implies $a^1= ax^1$ since $\cala_\rho$ is a domain. Thus one has $a^k=ax^k$, $\forall k \in \{1,\dots, 7\}$. $\square$

This implies that $\cala_\rho$ is Koszul and therefore is regular and is in fact a 3-Calabi-Yau domain.\\
\\

The matrix $A_{\beta}(x)$ is given by (\ref{AB7}) and the relations of the algebra $\cala_\beta$ read

\begin{equation}
[x^3,x^5]     +      [x^7,x^2]     =      0        \label{par1}
\end{equation}
\begin{equation}
[x^4,x^5]     +     [x^1,x^7]     =     0      \label{par2}
\end{equation}
\begin{equation}
[x^5,x^1]     +     [x^6,x^4]     =     0        \label{par3}
\end{equation}
\begin{equation}
[x^5,x^2]     +     [x^3,x^6]     =     0       \label{par4}
\end{equation}
\begin{equation}
[x^1,x^3]     +     [x^2,x^4]     =     0       \label{par5}
\end{equation}
\begin{equation}
[x^4,x^3]        =     0            \label{par6}
\end{equation}
\begin{equation}
[x^2,x^1]        =     0           \label{par7}
\end{equation}
where it is understood that the relations are valid in the algebra $\cala_\beta$ for the generator $x^\mu,\mu\in\{1,\dots,7\}$. The relations (\ref{par5}),  (\ref{par6}),  (\ref{par7}) between $x^1,x^2, x^3, x^4$ define the quadratic sub-algebra $\calb$ of $\cala_\beta$ which is Koszul of global dimension 2 and the following lemma implies in particular the injectivity of the $\calb$-module homomorphism
\[
\calb^3\stackrel{B(x)}{\longrightarrow}\calb^4
\]
where $B(x)$ is the $3\times 4$ matrix defined by \ref{B7} in the last section.
\begin{lemma}\label{Binv}
One has the following identity
\begin{equation}\label{invb}
\left(
\begin{array}{cccc}
-x^3 & -x^4 &x^1 & x^2\\
0&0&x^4&-x^3\\
x^2 &-x^1 &0 & 0
\end{array}
\right)
\left(
\begin{array}{ccc}
-x^1 & 0 & 2x^4\\
-x^2 & 0 & -2x^3\\
x^3 & 2x^2 & 0\\
x^4 & -2x^1 &0
\end{array}
\right) =
\left(
\begin{array}{ccc}
1&0&0\\
0&1&0\\
0&0&1
\end{array}
\right) u
\end{equation}
where $u\in \calb\subset \cala_\beta$ is given by
\begin{equation}\label{u}
u=2(x^1x^3+x^2x^4)=2(x^3x^1+x^4x^2)=x^1x^3+x^2x^4+x^3x^1+x^4x^2
\end{equation}
the different equalities follows from (\ref{par5}).
\end{lemma}

\noindent \underbar{Proof}. This is a direct consequence of the relations (\ref{par5}), (\ref{par6}), (\ref{par7}) which characterize $\calb\subset \cala_\beta$. $\square$\\

Let $(a^1, a^2, a^3, a^4, a^5, a^6, a^7)\in \cala^7_\beta$ be in the kernel of the $\cala_\beta$-module homomorphism

\begin{equation}\label{Ker beta}
\cala^7_\beta \stackrel{A_\beta(x)}{\longrightarrow} \cala^7_\beta
\end{equation}
that means that one has
\begin{equation}\label{ker beta}
(a^1,a^2,a^3,a^4,a^5,a^6,a^7)A_\beta(x)=0
\end{equation}
for the $a^\mu\in\cala_\beta$, $\mu\in\{1,\dots,7\}$. Let us cut the above relation (\ref{ker beta}) in two pieces as
\begin{equation}\label{kbeta1}
\begin{array}{l}
(a^1,a^2,a^3,a^4)\left(
\begin{array}{cccc}
0&x^7&-x^5&0\\
-x^7&0&0&-x^5\\
x^5&0&0&x^6\\
0&x^5&-x^6&0
\end{array}
\right)+ \\
\\
+ (a^5,a^6,a^7)\left(
\begin{array}{cccc}
-x^3 &-x^4 & x^1 & x^2\\
0&0&x^4 &-x^3\\
x^2 & -x^2 & 0&0
\end{array} \right)=0
\end{array}
\end{equation}
and 
\begin{equation} \label{kbeta2}
(a^1,a^2,a^3, a^4)\left(
\begin{array}{ccc}
x^3 & 0 & -x^2\\
x^4 & 0 & x^1\\
-x^1 & -x^4 &0\\
-x^2 & x^3 & 0
\end{array} \right)= 0
\end{equation}
which separates $(a^1,a^2, a^3, a^4)$ and $(a^5, a^6,a^7)$. We now use Lemma \ref{Binv} to express $a^5, a^6, a^7$ in terms of $a^1, a^2, a^3, a^4$. By applying the $4\times 3$-matrix
\[
\left(
\begin{array}{ccc}
-x^1 &0&2x^4\\
-x^2 & 0 &-2x^3\\
x^3 & 2x^2 & 0\\
x^4 & -2x^1 & 0
\end{array}
\right)
\]
on both terms of (\ref{kbeta1}), one obtains
\begin{equation}\label{a-a}
(a^1,a^2, a^3, a^4) C = (a^5, a^6, a^7)u
\end{equation}
where the $4\times 3$-matrix $C$ which coefficients in $\cala_\beta$ is given by
\begin{equation}\label{mC}
C= \left(
\begin{array}{ccc}
-(x^7x^2 + x^5x^3) & -2x^5x^2 & -2x^7x^3\\
x^7x^1 - x^5 x^4 & 2x^5x^1 & -2x^7x^4\\
-x^5x^1+x^6x^4 & -2x^6x^1 & 2x^5x^4\\
-(x^5x^2+x^6x^3)& -2x^6x^2 & -2x^5x^3
\end{array}
\right)
\end{equation}

 It is obvious that a solution of (\ref{kbeta2}) is $(a^1a^2a^3a^4)=a(x^1,x^2,x^3, x^4)$ with $a\in\cala_\beta$. It follows then from (\ref{a-a}), (\ref{mC}) that one has the implication
 \begin{equation}\label{Ia}
(a^1,a^2,a^3,a^4)=a(x^1, x^2,x^3, x^4)\Rightarrow (a^5,a^6,a^7)=a(x^5,x^6, x^7)
\end{equation}
as easily verified by using the relations (\ref{par1}), (\ref{par2}), (\ref{par3}),(\ref{par4}).

\begin{lemma}\label{bigr}
The algebra $\cala_\beta$ is bigraded by giving bidegree (1,0) to $x^1,x^2,x^3, x^4$ and bidegree (0,1) to $x^5, x^6, x^7$.
\end{lemma}

\noindent\underbar{Proof}. This is clear since the relations are homogeneous in the bidegree either of bidegree (2,0) or of bidegree (1,1). $\square$\\

In the following we shall refer to the degree in $(x^5,x^6,x^7)$ as the {\sl second degree} and our proof of the regularity of the quadratic algebra $\cala_\beta$ will be based on the induction with respect to this second degree $p\in \mathbb N$.  More precisely we shall prove by induction on the second degree $p\in \mathbb N$ the following statement.
\begin{proposition}\label{P4}
Let $a^1,a^2,a^3, a^4$ be 4 elements of $\cala_\beta$. Then $(a^1,a^2,a^3,a^4)$ satisfies Relation (\ref{kbeta2}) if 	and only if one has
\[
(a^1,a^2,a^3, a^4)=a(x^1,x^2, x^3, x^4)
\]
for some $a\in\cala_\beta$. 
\end{proposition}
The following lemma is the step 0 of the induction.
\begin{lemma}\label{s0}
Assume that $a^1, a^2, a^3, a^4\in \cala_\beta$ are of  second degree 0. Then $(a^1,a^2,a^3,a^4)$ satisfies Relation (\ref{kbeta2}) if and only if one has 
\[
(a^1,a^2,a^3,a^4)=a(x^1,x^2,x^3,x^4)
\]
for some $a\in\cala_\beta$ of second degree 0.
\end{lemma}

\noindent \underbar{Proof}. To assume that $a^k$ is of second degree 0 is the same as to assume $a^k\in\calb\subset \cala_\beta$. Since $\calb$ is Koszul of global dimension 2, one has the exact sequence
\[
0\longrightarrow \calb^3 \stackrel{B(x)}{\longrightarrow}\calb^4\stackrel{x}{\longrightarrow}\calb \longrightarrow \mathbb K\longrightarrow 0
\]
from which the result will follow. Indeed Relation (\ref{kbeta2}) reads
\[
\begin{array}{l}
(-a^3,-a^4,a^1, a^2)x=0\\
(0,0,a^4,-a^3)x=0\\
(a^2,-a^1,0,0)x=0
\end{array}
\]
which in $\calb$ is equivalent to
\[
\begin{array}{l}
(-a^3,-a^4,a^1,a^2)=(u^1,u^2,u^3)B(x)\\
(0,0,a^4,-a^3)=(v^1,v^2,v^3) B(x)\\
(a^2,-a^1,0,0)=(w^1,w^2,w^3)B(x)
\end{array}
\]
in view of the above exact sequence. From the second relation and the fact that $\calb$ is a domain follow that $(a^3,a^4)=r(x^3,x^4)$ and from the third relation that $(a^1,a^2)=s(x^1,x^2)$ with $r,s\in\calb$. Finally by using the first relation follow $r=s(=a)$. $\square$\\

Let us come back to the general case.
The equation (\ref{kbeta2}) is equivalent the equations
\begin{equation}\label{ax1}
-a^3x^1-a^4x^2+a^1x^3+a^2x^4=0
\end{equation}
\begin{equation}\label{ax2}
a^4x^3 -a^3x^4=0
\end{equation}
\begin{equation}\label{ax3}
a^2x^1-a^1x^2=0
\end{equation}
Equation (\ref{ax1}) is equivalent to
\begin{equation}\label{au}
(-a^3,-a^4,a^1,a^2,0,0,0)=(b^1,\dots,b^7)A_\beta(x)
\end{equation}
for some elements $b^1,\dots,b^7$ of $\cala_\beta$. It follows that one has
\begin{equation}
\begin{array}{l}\label{ab}
a^1=b^5x^1+b^6x^4-b^1x^5-b^4x^6\\
a^2=b^5x^2-b^6x^3-b^2x^5+b^3x^6\\
a^3=b^5x^3-b^7x^2-b^3x^5+b^2x^7\\
a^4=b^5x^4+b^7x^1-b^4x^5-b^1x^7
\end{array}
\end{equation}
and that 
\begin{equation}\label{bx}
(b^1,b^2,b^3,b^4)\left(
\begin{array}{ccc}
x^3 & 0 & -x^2\\
x^4 & 0 & x^1\\
-x^1 & -x^4 & 0\\
-x^2 & x^3 & 0
\end{array}
\right)
=0
\end{equation}

We now assume (induction hypothesis) that for any $a^i\in \cala_\beta$ $i\in\{1,2,3,4\}$ of second degree $q\leq p$ the relation (\ref{kbeta2})
implies 
\[
(a^1,a^2, a^3, a^4)=a(x^1,x^2, x^3, x^4)
\]
 for some $a\in\cala_\beta$, and let us assume that the $a^i$ in (\ref{ab}) are of second degree $q=p+1$. Then it follows that the $b^i$ are of second degree $q=p$ for $i\in \{1,2,3,4\}$ and therefore that $b^i=bx^i$ for some $b\in \cala_\beta$ for $i\in \{1,2,3,4\}$ in view of (\ref{bx}) and of the induction hypothesis. Then (\ref{ab}) implies
\begin{equation}\label{a1u}
a^1=(b^5-bx^5)x^1+(b^6-bx^6)x^4
\end{equation}
\begin{equation}\label{a2u}
a^2=(b^5-bx^5)x^2-(b^6-bx^6)x^3
\end{equation}
\begin{equation}\label{a3u}
a^3=(b^5-bx^5)x^3-(b^7-bx^7)x^2
\end{equation}
\begin{equation}\label{a4u}
a^4=(b^5-bx^5)x^4+(b^7-bx^7)x^1
\end{equation}
where we have used the relations (\ref{par3}), (\ref{par4}), (\ref{par1}) and (\ref{par2}). That is $a^1=ax^1+rx^4,a^2=ax^2-rx^3$, $a^3=ax^3-sx^2$, and $a^4=ax^4+sx^1$ with $a,r,s\in \cala_\beta$. By insertion of these expressions again in (\ref{ax2}) and (\ref{ax3}) one obtains
\[
sx^1x^3+sx^2x^4 = s(x^1x^3+x^2x^4)=0
\]
\[
-rx^3x^1-rx^4x^2=-r(x^1x^3+x^2x^4)=0
\]
which implies $r=s=0$ since $\cala_\beta$ is a domain. So $a^i=ax^i$ for $i\in\{ 1,2,3,4\}$ which achieves the proof of $(a^1,a^2, a^3,a^4)=a(x^1,x^2,x^3, x^4)$ by induction on the second degree $p$. $\square$\\

From  the implication (\ref{Ia}) it follows finally that one has
\[
(a^1,\dots,a^7)=a(x^1,\dots,a^7)
\]
for some $a\in\cala_\beta$ whenever $(a^1,\dots,a^7)$ is in the kernel of $A_\beta(x)$. This implies that $\cala_\beta$ is Koszul and therefore is regular of global dimension 3 or equivalently here is 3-Calabi-Yau.\\

In summary we have proved the following result      
\begin{theorem}\label{T7}
The quadratic algebras associated with the 3-regular exterior 3-forms $\rho, \beta, \alpha^{(3)}, \alpha^{(3)'}, \alpha^{(3)''}$ in dimension 7 and with  the exterior 3-forms $\alpha^{(p)}$ in dimensions $2p+1$ for $p\in \mathbb N$ with $p\geq 1$ are regular and are therefore 3-Calabi-Yau domains.
\end{theorem}
     Using the fact that the 3-regular exterior 3-forms that we have considered contain representative elements in all orbits of 3-regular exterior 3-forms in $\mathbb K^n$ for $n\leq 7$ whenever $\mathbb K$ is algebraically closed, one has the following theorem.

\begin{theorem}\label{n7}
Assume that $\mathbb K$ is algebraically closed and let $\alpha$ be an exterior 3-regular 3-form on $\mathbb K^n$ with $n\leq 7$. Then the quadratic algebra
\[
\cala=\mathbb K\langle x^1,\dots,x^n\rangle/[\{ \partial_i\alpha\}]
\]
 is regular which implies that it is a 3-Calabi-Yau domain.
\end{theorem}

It is worth noticing here that for $n\leq 7$, we have shown that such exterior 3-regular 3-forms only exist in dimensions $n=3, 5$ and 7.\\

In the next paper of this series we shall investigate a sequence of quadratic algebras of the above type associated with the sequence of the simple complex Lie algebras.

\section{Exterior 3-forms in $\mathbb K^n$ with $n\geq 8$}\label{Cex}

A natural question is : Does Theorem \ref{n7} remain true for exterior 3-forms of rank $n>7\>$ ?  Or, in other words, is the quadratic algebra associated with a 3-regular exterior 3-form in $\mathbb K^n$ with $n\geq 8$  a Koszul algebra ?\\

It turns out that the answer is negative. Indeed in \cite{roo:2010}  it is pointed out that the case XII of rank 8 of the book of G.B. Gurevich \cite{gur:1964} namely the exterior 3-form
\[
\omega= \theta^2 \wedge \theta^1 \wedge \theta^7+\theta^2\wedge \theta^5\wedge \theta^3+\theta^2\wedge \theta^6\wedge \theta^4+\theta^8 \wedge\theta^1 \wedge \theta^4+\theta^7\wedge \theta^3\wedge \theta^4
\]
which is 3-regular leads to an associated quadratic algebra which is not Koszul and therefore not regular. Let us explain this fact.\\

First the 3-regularity of $\omega$ follows from a tedious but straightforward calculation.\\
The relations of the quadratic algebra $\cala_\omega=\cala$ associated with $\omega$ read
\begin{equation}\label{RA}
\left\{
\begin{array}{l}
\partial_1\omega : [x^2,x^7]+[x^8,x^4]=0\\
\partial_2\omega : [x^7,x^1]+[x^3,x^5]+[x^4,x^6]=0\\
\partial_3\omega : [x^7,x^4]+[x^5,x^2]=0\\
\partial_4\omega : [x^3,x^7]+[x^1,x^8]+[x^6,x^2]=0\\
\partial_5\omega : [x^2,x^3]=0\\
\partial_6\omega : [x^2,x^4]=0\\
\partial_7\omega : [x^1,x^2]+[x^4,x^3]=0\\
\partial_8\omega : [x^4,x^1]=0
\end{array}
\right.
\end{equation}

These relations show that $\cala$ is generated by two sub-algebras : a tensor algebra 
$\mathbb K\langle x^5,x^6,x^7,x^8\rangle = \calt$ generated by $x^5,x^6,x^7,x^8$ and a quadratic sub-algebra $\calb_\omega=\calb$ generated by $x^1,x^2,x^3,x^4$ with relations $\partial_5\omega,\partial_6\omega,\partial_7\omega,\partial_8\omega$ while the first four relations $\partial_1\omega,\partial_2\omega,\partial_3\omega,\partial_4\omega$ of $\cala$ give commutation relations between elements of $\calt$ and elements of $\calb$. The relations of $\cala$ read $\partial\omega=A(x)x$ as explained in Sections \ref{Sf} (formulas \ref{Ax} and \ref{relA} where $x$ is the column $\left(\begin{array}{c}
x^1\\
\vdots\\
x^8
\end{array}\right)$
and where the antisymmetric $8\times 8$-matrix $A(x)=A_\omega(x)$ is given by
\begin{equation}\label{Aom}
A(x)=\left[
\begin{array}{cccccccc}
0 & x^7 & 0 & -x^8 & 0 & 0 & -x^2 & x^4\\
-x^7 & 0 & x^5 & x^6 & -x^3 & -x^4 & x^1 & 0\\
0 & -x^5 & 0 & -x^7 & x^2 & 0 & x^4 & 0\\
x^8 & -x^6 & x^7 & 0 & 0 & x^2 & -x^3 & -x^1\\
0& x^3 & -x^2 & 0 & 0 & 0 & 0 & 0\\
0& x^4 & 0 & -x^2 & 0 & 0 & 0 & 0\\
x^2 & -x^1 &-x^4 & x^3 & 0 & 0 & 0 & 0\\
-x^4 & 0 &0 & x^1 & 0 & 0 & 0 & 0
\end{array}
\right]
\end{equation}
which reads in $4\times 4$-matrices blocks
\begin{equation}\label{Aom2}
A(x)=\left(
\begin{array}{cc}
A_0(x_{II}) &-B(x_I)^t\\
\\
B(x_I) & 0
\end{array}
\right)
\end{equation}
where $x_I=(x^1,x^2,x^3,x^4)^t$ and $x_{II}=(x^5,x^6,x^7,x^8)^t$.\\
The 3-regularity of $\omega$ implies \cite{mdv:2007} that the augmented Koszul complex of $\cala$ reads
\[
0\rightarrow \cala\otimes \omega \stackrel{d}{\rightarrow} \cala\otimes R\stackrel{d}{\rightarrow} \cala\otimes E \stackrel{d}{\rightarrow} \cala\rightarrow \mathbb K \rightarrow 0
\]
where $E=\cala_1$ is the space generated by the $x^1,\dots, x^8$ and where $R\subset E\otimes E$ is the space of relations of $\cala$ generated by the $\partial_1\omega,\dots,\partial_8 \omega$ while $\omega \in \wedge^3E\subset E\otimes E\otimes E$ is identified as an element of $E^{\otimes^3}$. With other notations the above Koszul complex of $\cala$ identifies to the sequence (see in Section \ref{Sf})
\[0\rightarrow \cala\stackrel{x^t}{\rightarrow} \cala^8 \stackrel{A(x)}{\rightarrow}\cala^8\stackrel{x}{\rightarrow} \cala\rightarrow \mathbb K\rightarrow 0
\]
where the exactness of
\[
\cala^8\stackrel{A(x)}{\rightarrow}\cala^8\stackrel{x}{\rightarrow}\cala\rightarrow \mathbb K\rightarrow 0
\]
is automatic since it is equivalent to the presentation of $\cala$ by generators and relations while the exactness of
\[
0\rightarrow \cala\stackrel{x^t}{\rightarrow} \cala^8
\]
follows from the fact that $\cala$ is a domain. Thus it remains to analyse the small complex
\begin{equation}\label{SC8}
\cala\stackrel{x^t}{\rightarrow} \cala^8\stackrel{A(x)}{\rightarrow} \cala^8
\end{equation}
in other words to describe the kernel $\ker(A(x))$ of $A(x)$. By construction one has
\begin{equation}\label{K1}
a_1(x^1,\dots, x^8)\in \ker (A(x)),\> \>  \forall a_1\in \cala
\end{equation}
since (\ref{SC8}) is a complex, however there is another injective mapping of $\cala$ into $\ker(A(x))$ which is quadratic in the $x^i$ defined by
\begin{equation}\label{K2}
a_2(j^1(x),\dots,j^8(x))\in \ker(A(x)),\>\>  \forall a_2\in \cala
\end{equation}
where
\begin{equation}\label{Ku0}
j^1=j^2=j^3=j^4=0
\end{equation}
\begin{equation}\label{Ku1}
j^5=u^1(x_I)=-(x^4)^2
\end{equation}
\begin{equation}\label{Ku2}
j^6=u^2(x_I)=x^2x^1+x^4 x^3+[x^2,x^1] \>\> (\simeq (x^2x^1+x^4x^3+[x^4,x^3]))
\end{equation}
\begin{equation}\label{Ku3}
j^7=u^3(x_I)=x^2x^4( \simeq x^4x^2)
\end{equation}
\begin{equation}\label{Ku4}
j^8=u^4(x_I)=(x^2)^2
\end{equation}
(the equivalences are modulo the relations $\partial_i\omega$ for $i\in \{5 ,6,7,8\}$). It is easy to verify that that (\ref{K1}) and (\ref{K2}) belong to 
$\ker(A(x))$ and one can show that they generate $\ker(A(x))$. This implies in particular that the sequence (\ref{SC8}) is not exact and therefore that $\cala$ is not a Koszul algebra and so cannot be regular. Notice that nevertheless one has the exact sequence
\begin{equation}\label{Sd3}
0\longrightarrow \cala^2\stackrel{(x^t,j)}{\longrightarrow} \cala^8\stackrel{A(x)}{\longrightarrow} \cala^8 \stackrel{x}{\longrightarrow}\cala \longrightarrow \mathbb K\longrightarrow 0
\end{equation}
which is a minimal projective resolution of the $\cala$-module $\mathbb K$. This implies  \cite{car:1958} that the global dimension of $\cala$ is 3, $\gldim (\cala)=3$.\\

The origin of the above facts is that the global dimension of the sub-algebra $\calb$ is 3, $\gldim(\calb)=3$ and that one has the minimal free resolution of $\mathbb K$ as $\calb$-module
\begin{equation}\label{SB}
0\longrightarrow \calb \stackrel{u}{\longrightarrow}\calb^4\stackrel{B(x_I)}{\longrightarrow} \calb^4 \stackrel{x_I}{\longrightarrow}\calb \longrightarrow \mathbb K \longrightarrow 0
\end{equation}
since the mapping $u:\calb\rightarrow \calb^4$ given by 
\begin{equation}\label{KB}
b\mapsto b(u^1(x_I), u^2(x_I), u^3(x_I), u^4(x_I))
\end{equation}
is injective and that $\im(u)=\ker(B(x_I))$ as easily verified.\\

By comparison of the above resolutions with the usual form of the minimal projective resolutions, (i.e. here minimal free resolutions since we are in the graded case \cite{car:1958})
\[
\dots \rightarrow \cala\otimes \Tor_n^\cala (\mathbb K,\mathbb K) \stackrel{d_n}{\rightarrow}\dots \stackrel{d_1}{\rightarrow} \cala\otimes \Tor^\cala_0(\mathbb K, \mathbb K)\rightarrow \mathbb K\rightarrow 0
\]
one gets that
\[
\begin{array}{lll}
\Tor_3^\cala(\mathbb K, \mathbb K) & = & \Tor^\cala_{3,3}(\mathbb K, \mathbb K)\oplus \Tor^\cala_{3,4}(\mathbb K, \mathbb K)\\
\\
\Tor_3^\calb(\mathbb K, \mathbb K) & = & \Tor^\calb_{3,4}(\mathbb K, \mathbb K)
\end{array}
\]
with
\[
\dim(\Tor^\cala_{3,3}(\mathbb K, \mathbb K))=\dim(\Tor^\cala_{3,4}(\mathbb K, \mathbb K))=\dim(\Tor^\calb_{3,4}(\mathbb K, \mathbb K))=1
\]
and
\[
\Tor^\cala_n(\mathbb K, \mathbb K)=\Tor^\calb_n(\mathbb K,\mathbb K)=0
\]
for $n\geq 4$. It follows that the Poincaré double series 
\begin{equation}\label{Pg}
P_\cala(r,s)=\sum_{k,\ell}\dim(\Tor^\cala_{k,\ell}(\mathbb K,\mathbb K)) r^ks^\ell
\end{equation}
of $\cala$ and $\calb$ read
\begin{equation}\label{PA}
P_\cala(r,s)= 1 + 8 rs+8r^2s^2+r^3s^3+r^3s^4
\end{equation}
\begin{equation}\label{PB}
P_\calb(r,s)= 1 + 4 rs + 4r^2s^2 + r^3s^4
\end{equation}
from which one obtains the Hilbert series
\begin{equation}\label{Hg}
H_\cala(t)=\sum_n \dim(\cala_n) t^n
\end{equation}
of $\cala$ and $\calb$
\begin{equation}\label{HA}
H_\cala(t)=1/1-8 t + 8 t^2- t^3-t^4
\end{equation}
which is the result of \cite{roo:2010}, and 
\begin{equation}\label{HB}
H_\calb(t)=1/1-4 t+4 t^2- t^4
\end{equation}
by using the general formula
\begin{equation}\label{HPg}
H_\cala(t) P_\cala(-1,t)=1
\end{equation}
and the above results.

\section{Conclusion and further prospects}

It is worth noticing here that for the 3-linear forms which are 3-regular there is already an example in dimension 3 for which the associated quadratic algebra is not regular \cite{mdv:2010} which shows that the conjecture of \cite{mdv:2007} is wrong. However this counterexample is very asymmetric from the point of view of the group of permutations between the coordinates. This is why  it was interesting to investigate the case of the 3-regular exterior 3-forms for which the complete antisymmetry is a maximal symmetry for the group of permutations between the coordinates. There we have shown that the first counter-example occurs only in dimension 8.\\

Nevertheless it is interesting to consider natural families of 3-regular exterior 3-forms as  the family of 3-regular $\alpha^{(p)}\in \wedge^3\mathbb K^{2p+1}$ investigated in Section 3 which leads to associated regular quadratic algebras (see Proposition 5 in Section 3). A very interesting such family is the family of canonical 3-forms of simple Lie algebras. Let us explain this point.\\

Let $\fracg$ be a finite-dimensional Lie-algebra and let $\omega$ be the 3-linear form on $\fracg$ defined by
\begin{equation}\label{3f}
\omega(X,Y,Z)=\Tr([\ad(X),\ad(Y)]\ad(Z))
\end{equation}
for  $X,Y,Z\in \fracg$, where $\Tr(.)$ denotes the trace of the linear endomorphisms of $\fracg$. By definition $\omega(X,Y,Z)$ is antisymmetric in $X$ and $Y$, furthermore the symmetries of the trace implies
\begin{equation}\label{cyclf}
\omega(X,Y,Z)=\omega(Z,X,Y)
\end{equation}
for any $X,Y,Z\in \fracg$ as well as the $\ad$-invariance of $\omega$. Thus $\omega$ is an invariant exterior 3-form on $\fracg$ $(\omega\in \wedge^3_I \fracg^\ast)$ which will be referred to as {\sl the canonical 3-form of $\fracg$} (the above  terminology is the same as the one of \cite{kab:2009} while  in \cite{hvl:2013} for instance, $\omega$ is referred to as the {\sl Cartan 3-form of $\fracg$}). This invariant exterior 3-form is closely related to the structure constants of $\fracg$. Indeed let $(E_k)$ be an arbitrary basis of $\fracg$, the relations of $\fracg$ read in this basis
\begin{equation}\label{rel g}
[E_k, E_\ell]=C_{k\ell}^{\>\>\>\>r}E_r
\end{equation}
where the $C_{k\ell}^{\>\>\>\>r}\in \mathbb K$ are the corresponding structure constants of $\fracg$. Then the components $\omega_{k\ell m}=\omega(E_k,E_\ell,E_m)$ of $\omega$ read

\begin{equation}\label{alf str}
\omega_{k\ell m}=C_{k\ell}^{\>\>\>\>r} K_{rm}
\end{equation}
where $K_{rm}=K(E_r,E_m)$ are the components of the {\sl the Killing form} of $\fracg$ defined by
\begin{equation}\label{Kill}
K(X,Y)=\Tr (\ad(X)\ad(Y))
\end{equation}
for any $X,Y\in \fracg$.\\

In the sequel we assume that $\mathbb K= \mathbb C$ so that the Lie algebra $\fracg$ is a complex Lie algebra.

\begin{lemma}\label{SSLa}
The canonical 3-form $\omega$ of $\fracg$ is nondegenerate if and only if $\fracg$ is a semisimple Lie algebra.
\end{lemma}
\noindent {\bf Proof}. This is a consequence of the following connection between $\omega$ and the Killing form of $\fracg$
\begin{equation}\label{aK}
\omega(X,Y,Z)=K(X,[Y,Z])
\end{equation}
for any $X,Y,Z\in \fracg$. Indeed if $i_X\omega=0$ then $K(X,[YZ])=0$ for  any $Y,Z\in\fracg$ so if $i_X\omega=0$ implies $X=0$, then $K(X,[Y,Z])=0$ for any $Y,Z\in \fracg$ implies $X=0$. Therefore if $\omega$ is nondegenerate then a fortiori the Killing form is nondegenerated so that $\fracg$ is semi-simple. Conversely if $\fracg$ is semi-simple $\fracg=[\fracg,\fracg]$ so the relation (\ref{aK}) implies that the nondegeneracy of $\omega$ is equivalent to the one of $K$. $\square$\\

The content of Lemma \ref{SSLa} above is the same as the one of Lemma 2.1 of \cite{hvl:2013}.\\

As already pointed out, $\omega$ nondegenerate is equivalent to $\omega$ preregular for an exterior 3-form.\\

\begin{theorem}\label{ThRg}
The canonical 3-form $\omega$ of $\fracg$ is 3-regular if and only if $\fracg$ is a simple Lie algebra.
\end{theorem}
\noindent {\bf Proof}. The Lie algebra $\fracg$ is simple if and only if its adjoint representation is irreducible. This irreducibility is equivalent in the present context to the fact that the commutant $\ad(\fracg)'$ of $\ad(\fracg)$ consists of multiples of the identity mapping of $\fracg$ onto itself, that is to the condition

\begin{equation}\label{adirr}
\begin{array}{l}
[L, \ad X]=0,\>\>\> \forall X\in \fracg\\
\text{implies}\\
L=\lambda I \> \> \text{with}\>\> \lambda\in \mathbb C
\end{array}
\end{equation}
for any linear endomorphism $L$ of $\fracg$ where I denotes the identity mapping of $\fracg$.\\

One may assume that $\fracg$ is semisimple since both $\omega$ 3-regular or $\fracg$ simple implies the semisimplicity of $\fracg$. The Killing form $K$ of $\fracg$ is then nondegenerate. Now, 
\[
\omega(LX,Y,Z)-\omega(X,MY,Z)=0, \>\>\> \forall X,Y,Z\in \fracg
\]
is equivalent to
\[
K(Z,[LX,Y]-[X,MY])=0,\>\>\> \forall X,Y,Z\in \fracg
\]
in view of the relation (\ref{aK}) and finally is equivalent to 
\[
[LX,Y]-[X,MY]=0, \>\>\> \forall X,Y\in \fracg
\]
since $K$ is nondegenerate. Thus the assumption that $\omega$ is 3-regular is equivalent to the condition
\begin{equation}\label{eq3r}
\begin{array}{l}
[LX,Y]-[X,MY]=0,\>\>\> \forall X,Y\in \fracg\\
\text{implies}\\
L=M=\lambda I\> \text{with}\> \lambda\in \mathbb C
\end{array}
\end{equation}
for the linear endomorphisms $L$ and $M$ of $\fracg$, where $I$ denotes the identity mapping of $\fracg$.\\

Let us first prove that the 3-regularity of $\omega$ implies the simplicity of $\fracg$, that is that one has the implication (\ref{eq3r})$\Rightarrow$ (\ref{adirr}). By setting $M=L$ in condition (\ref{eq3r}) one has that
\[
\ad (X) LY +\ad (Y)LX=0
\]
implies $L=\lambda I$. Now if $L\in \ad(\fracg)'$ one has
\[
\ad (X)LY + \ad (Y) L(X) = L([X,Y] +[Y,X])=0
\]
and therefore $L=\lambda I$. This means that $\fracg$ is simple whenever its canonical 3-form $\omega$ is 3-regular.\\

Let us show that conversely the simplicity of $\fracg$ implies the 3-regularity of its canonical 3-form $\omega$.\\

So let $\fracg$ be a simple $n$-dimensional complex Lie algebra and let the $E_k$ $(k\in \{1,\dots, n\})$ be a basis of $\fracg$ which is orthonormal for the Killing form i.e. such that
\begin{equation}\label{ortn}
K(E_k, E_\ell) = \Tr (\ad(E_k)\ad(E_\ell))=\delta_{k\ell}
\end{equation} 
for  $k,\ell\in \{1,\dots,n\}$. The relations of $\fracg$ read then
\begin{equation}\label{reLie}
[E_k,E_\ell]=\sum_m \omega_{kl\ell m} E_m
\end{equation}
where the $\omega_{k\ell m}$ are the components of the canonical 3-form $\omega$ in the basis $(E_k)$. The matrix components of $\ad(E_k)$ are given by
\begin{equation}\label{matad}
\ad(E_k)^m_\ell = \omega_{k\ell m}
\end{equation}
$\forall k,\ell, m\in \{1,\dots, n\}$ and the orthonormality of the $E_k$ reads in terms of the $\omega_{k\ell m}$
\begin{equation}\label{alf ortn}
\sum_{ij} \omega_{ijk} \omega_{ij\ell} = \delta_{k\ell}
\end{equation}
$\forall k,\ell\in \{1,\dots,n\}$ while the Jacobi identity reads
\begin{equation}\label{Jac}
\sum_{m}(\omega_{ijm} \omega_{k\ell m}+ \omega_{kim} \omega_{j\ell m}+ \omega_{jkm} \omega_{i\ell m})=0
\end{equation}
$\forall i,j,k,\ell \in \{1,\dots,n\}$.\\

Let $L$ and $M$ be two linear endomorphisms of $\fracg$ satisfying
\begin{equation}\label{aLM}
\omega(LX,Y,Z)=\omega(X,MY,Z)\>\>\>\> \forall X,Y,Z \in \fracg
\end{equation}
which reads in components
\begin{equation}\label{aLMc}
L^r_i \omega_{rjk}=M^s_j \omega_{isk}
\end{equation}
$\forall i,j,k$ and let us show that this implies that $L=M=\lambda I$. The antisymmetry of $\omega_{ijk}$ implies
\begin{equation}\label{aaLMc}
L^r_i\omega_{rji} =0
\end{equation}
while (\ref{alf ortn}) implies
\[
L^r_i=\sum_{jk} M^s_j \omega_{isk} \omega_{rjk} \Rightarrow \Tr(L) = \Tr(M)
\]
which leads by using (\ref{Jac}) and (\ref{aaLMc}à) to
\[
L^r_i = \sum_{jk} M^s_j \omega_{srk}\omega_{ijk} = L^t_r \omega_{tjk} \omega_{ijk}=L^i_r
\]
for $r,i\in \{1,\dots, n\}$. Similarily one gets $M^s_j=M^j_s$ for $s,j\in \{1,\dots,n\}$.  This implies
\[
\ad(E_k)L=M\ad (E_k)
\]
and by transposition
\[
\ad(E_k)M=L\ad(E_k)\>\>\>\> \forall k
\]
so one has $[\ad(E_k), L+M]=0$ which by irreducibility of the adjoint representation (simplicity of $\fracg$) implies $L+M=2\lambda I$ for some $\lambda\in \mathbb C$. Then
\[
\ad(E_k)(L-\lambda I)=-(L-\lambda I)\ad (E_k)
\]
which implies $[\ad(E_k),\ad(E_\ell)] (L-\lambda I)=0$
so finally (by simplicity)
\[
L=M=\lambda I
\]
which means that the canonical 3-form $\omega$ is 3-regular. $\square$\\

Thus the family of canonical 3-forms of simple Lie algebras is a family of 3-regular exterior 3-forms. It is natural to investigate the regularity of the associated quadratic algebras. For the lowest dimensional simple Lie algebra $\pfraca_1=\fracsl(2)$ which is of dimension 3 it is straightforward to show that the associated quadratic algebra is the polynomial algebra $\mathbb C[x^1, x^2, x^3]$ which is regular of global dimension 3. For the next simple Lie algebra $\pfraca_2=\fracsl(3)$ which is of dimension 8, one can show with a computer that the associated quadratic algebra is Koszul of global dimension 3 but an analytic proof is not straightforward. To go further on one must use the root space decomposition since it is clear that these quadratic algebras are in fact associated with  the irreducible root systems. This work is in progress.

\bigskip
\bigskip
\bigskip
\bigskip

\noindent {\bf Acknowledgements}

\medskip

Work supported by the project MTM2017-86987-P  from MICIN (Spain) and the group FQM 211 from Junta de Andalucía.\\
MDV thanks Roland Berger for discussion.

\newpage 

\bibliographystyle{plain}
\bibliography{BibMich,BibExt}

\end{document}